\documentclass[a4paper,12pt]{article}
\input{isolatin1.sty}

\newcommand{\relbd}{{\rm relbd}}
\newcommand{\conv}{{\rm conv\,}}
\newcommand{\aff}{{\rm aff\,}}

\newcommand{\R}{{\Bbb R}}

\newcommand{\Ol}[1]{\overline{#1}}

\newcommand{\qed} {\hfill $\Box$ \protect\bigskip \par}

\usepackage{epsf}
\usepackage{epsfig}
\usepackage{latexsym}
\usepackage{amsmath,amsfonts,amssymb}
\newtheorem{lemma}{Lemma}[section] 

\newtheorem{satz}[lemma]{Theorem}

\newtheorem{prop}[lemma]{Proposition}
 
\newtheorem{algo}[lemma]{Algorithm}

\newtheorem{defi}[lemma]{Definition}
\begin{document}
\title{Relaxation, New Combinatorial and Polynomial Algorithms for the Linear Feasibility Problem}
\author{Ulrich Betke\thanks
{Fachbereich Mathematik, Universität Siegen, D--57068 Siegen, Germany, e--mail: betke@mathematik.uni--siegen.de}}
\maketitle
\newcommand{\defe}{{\rm def\,}}
\newcommand{\co}{{\rm co\,}}

\begin{abstract}
We consider the homogenized linear feasibility problem, to find an $x$ on  the unit sphere, satisfying $n$ linear inequalities $a_i^Tx\ge 0$. To solve this problem we consider the centers of the insphere of spherical simplices, whose facets are determined by a subset of the constraints. As a result we find a new combinatorial algorithm for the linear feasibility problem. If we allow rescaling this algorithm becomes polynomial. We point out that the algorithm solves as well the more general convex feasibility problem. Moreover numerical experiments show that the algorithm could be of practical interest.


\end{abstract}

\section{Introduction}
A central problem in optimization is the {\em linear feasibility problem} F. An instance $A$ of F is to find for $n$  vectors $a_i\in \R^d$, $i=1,\dots,n$ and $n$ real numbers $b_i$, $i=1,\dots,n$
\begin{equation}\label{L}
x\in \R^d \text{ such that } a_i^Tx\ge b_i,\quad i=1,\dots,n
\end{equation}
or to show that such an $x$ does not exist.

For an instance $A$ of F the set $P_A$ of {\em feasible points}
\begin{equation}\label{P_L}
P_A=\{x \in \R^d \mid a_i^T x \ge b_i,\;i=1,\dots,n\}
\end{equation}
is a polyhedron.

The main application of F is the {\em linear optimization problem} to minimize a linear function over $P_A$. It is well known how to solve an instance of the linear optimization problem by solving of one or more instances of F. Of the vast literature on linear optimization we mention the books \cite{GLS}, \cite{PS}, \cite{S}, which provide the background for this paper. 

The linear feasibility problem is a particular case of the convex feasibility problem. Here we want to find a point in a convex set $K$. Clearly the problem depends on the representation of $K$. A very general way of representing $K$ is by a separation oracle \cite{GLS}. For every $x\in \R^d$ this oracle either confirms that $x\in K$ or gives a hyperplane separating $x$ from $K$.  A  survey on algorithms for the convex feasibility problem  is \cite{BB}.

There are numerous algorithms for F. Well known algorithms  are the {\em phase I} of the {\em simplex algorithm}, {\sc Khachian}'{\em s algorithm} \cite{Kh} and {\sc Karmar\-kar}'{\em s algorithm} \cite{Ka}, for a fairly recent survey of {\sc Karmarkar}'s algorithm see e.g. \cite{PW}. Here  {\sc Khachian}' s algorithm and {\sc Karmarkar}'s algorithm are distinguished by the fact that they are polynomial. The simplex algorithm is distinguished by the fact that it is combinatorial, i.e. it's state is completely determined by a subset of the constraints.

For our purposes {\sc Khachiyan}'s algorithm and an algorithm given by  {\sc Agmon}, {\sc Motzkin} and {\sc Schoenberg} \cite{A}, \cite{MS} are of particular interest as they are {\em relaxation algorithms} for F which may be described by the following simple prototype

\begin{algo}\label{relax}
\begin{enumerate}  
\item
Choose an arbitrary $x^0$.
\item
Choose a ``simple'' convex set $C_k$ depending on $x^{k-1}$ and other data obtained in the course of the algorithm, such that $P_A\subset C_k$ and a suitable point $x^k\in C_k$. If $x^k\in P_A$ terminate.
\item Set $k:=k+1$ and goto step (2.).
\end{enumerate}
\end{algo}

In the case of {\sc Khachian}'s algorithm, $C_k$ is an ellipsoid and $x^k$ its center. The next ellipsoid is constructed in the following way. Choose a constraint $i$, such that $a_i^Tx^k < b_i$. Then the next ellipsoid $C_{k+1}$ is the smallest ellipsoid circumscribed to
$$
C_k \cap \{x \mid a_i^Tx\ge a_i^Tx^k\}.
$$

For an instance $A$ of F and $x\not \in P_A$ we say that $i\in \{1,\dots,n\}$ is a {\em most violated constraint}, if $a_i^Tx-b_i=\min\{a_j^Tx-b_j\mid j=1,\dots,n\}$. 
In  the {\sc Agmon--Motzkin--Schoenberg} algorithm a most violated constraint, $a_i^Tx^k\ge b_i$, say, for $x^{k}$ is chosen and $x^{k+1}$ is the projection of $x^{k}$ on the half-space $\{a_i^Tx\ge b_i\}$. This algorithm is a very special case of an algorithm in \cite{FZ}.

Furthermore the simplex algorithm for a problem in standard form may be viewed as a relaxation algorithm for the dual linear feasibility problem, where each iterate $x^k$ is determined by the requirement that it satisfies $d$ of the inequalities in (\ref{L}) with equality.

Further we remark that {\sc Khachiyan}'s algorithm solves the problem of constructing a point in a convex set given by a separation oracle. We shall return to this topic after the presentation of our algorithms.

We proceed as follows: In section 2 we present a relaxation algorithm for the homogenized form of F and determine its properties. In section 3 we use additional transformations of the problem to construct  a polynomial algorithm for F. Moreover we shortly discuss its relation to the convex feasibility problem. While in the previous sections we have concentrated on the geometric aspects of the algorithms, we discuss in section 4 the underlying linear algebra and present the results of some numerical experiments.

\section{Algorithms}
In a certain way we want to combine the algorithm of {\sc Agmon--Motzkin--Schoenberg} with {\sc Khachiyan}'s algorithm. This means that we want to use a circumscribed body of $P_A$ and a point, which we may consider as a projection on this body. Furthermore we want this body as closely related to $P_A$ as possible. A rather obvious choice is a simplex bounded by $d+1$ of the hyperplanes $\{a_i^Tx=b_i\}$ and the center of its inball, as this point may be considered as a simultaneous projection on all half-spaces. For several reasons this does not work very well in Euclidean space and thus we first use the well known process of homogenization to transform the problem into a spherical problem. 

To proceed we need some more notation. We write $B(c,r)=\{x\in \R^d \mid \|x-c\|\le r\}$ for the ball with center $c$ and radius $r$, $S(c,r)=\{x\in \R^d \mid \|c-x\|=r\}$ for the sphere bounding $B(c,r)$ and $S^{d-1}=S(0,1)$ for the unit sphere in $\R^d$.

Then we define the homogenization
\begin{equation}\label{Homogen}
\Phi :\R^d \to S^d, \quad \Phi(x_1,\dots,x_d)= (x_1,\dots,x_d,1)/\|(x_1,\dots,x_d,1)\|.
\end{equation}

For $x=(x_1,\dots,x_{d+1})\in \R^{d+1}$ we denote by $\Ol x=(x_1,\dots,x_d)$ its projection on $\R^d$. 
For an instance $A$ of F we have for the set of feasible points $P_A$
$$
\Phi(P_A) = \{x\in S^d \mid a_i^T \Ol x -b_i x_{d+1}\ge 0,\; i=1,\dots,n,\; x_{d+1}>0\}.
$$
Thus we obtain an algorithm for the linear feasibility problem F, if we have an algorithm for the {\em spherical feasibility problem} F', where an instance $A$ of F' is to find
\begin{equation}\label{P_S}
x\in S^{d} \text{ such that } a_i^T x \ge 0,\quad i=1,\dots, n
\end{equation}
or to show that such an $x$ does not exist. For an instance $A$ of F' we denote again by $P_{A}$ the set of its feasible points.

To study F' it is helpful to use the notion of {\em violation}. For an instance $A$ of F' and $x\in S^d$ we define the violation $v(x)$ by
$$
v(x)=\max\{0, \max\{-a_i^Tx \mid i=1,\dots,n\}\}.
$$
Thus a point is feasible for an instance $A$ of F, if and only if its violation is $0$.

For our algorithms it is somehow more appropriate to take a polar point of view, i.e. we work with the normals $a_i$ rather than with the half spheres $a_i^T x\ge 0$. Moreover 
 we  construct a sequence of points, which show that the instance is close to being infeasible. To be more precise we denote the origin of $\R^{d+1}$ by $O$, the {\em convex hull}, {\em affine hull} respectively, of a set $M\subset \R^{d+1}$ by  $\conv M$, $\aff M$ respectively, and define

\begin{defi} \label{pos_span}
Let $Q=\{x^i\in \R^{d+1} \mid i=1,\dots,k\}$ be a set of points. We say that $Q$ is {\em positively spanning} if $Q$ is affinely independent and $O\in \conv Q$.
\end{defi}

If  $k=d+1$ then $Q$ is a {\em positive basis}. Here we also consider the case that $\dim (\aff Q)<d+1$.

\begin{defi} Let $Q=\{x^i\in \R^{d+1} \mid i=1,\dots,k\}$ be a set of points. We say that $Q$ is {\em nearly positively spanning} if $Q$ is affinely independent and the orthogonal projection $O'$ of $O$ on $\aff Q$ is contained in $\conv Q$. The distance of $O'$ to the origin is called the {\em deficiency} of $Q$ and denoted by $\defe Q$.
\end{defi}

Thus a nearly positively spanning set is positively spanning if and only if its deficiency is 0. Moreover for a nearly positively spanning set the deficiency is equal to the distance of the origin to its convex hull.

The evident significance of positively spanning sets for F' is given by

\begin{prop} \label{Pos_Prop}
Let an instance $A$ of F' be given and $Q\subset \{a_1,\dots,a_n\}$ be a positively spanning set. Then every $x\in P_A$ satisfies $a_i^T x =0$ for all $a_i\in Q$. 
\end{prop} 

If the cardinality of $Q$ is $d+2$ then $A$ is infeasible. In the other case $P_A$ is contained in the plane $\{a_i^T x=0 \mid a_i\in Q\}$ and thus the dimension of the instance is reduced. In this case we may change the problem slightly in a way that $Q$ is transformed into a nearly positively spanning set with positive deficiency without changing the feasibility or infeasibility of the instance. This is quite analogous to cases of degeneration in {\sc Khachian}'s algorithm and will not be discussed further (but compare the remarks after Lemma \ref{Homo_vol}).

If the points of $Q$ are contained in $S^d$ then positively spanning sets can be nicely characterized.

\begin{defi} \label{touching}
Let $Q$ be a set of affinely independent points. The sphere $S(C,R)$ is said to be {\em touching} for $Q$, if $C\in \aff Q$ and $Q\subset S(C,R)$.
\end{defi}

Clearly every set of affinely independent points has a unique touching sphere. Moreover we observe that it is easy to compute a touching sphere for given $Q$. The relation of touching spheres and nearly positively spanning sets is given by

\begin{prop} \label{Um_Prop}
Let $Q \subset S^d$ be affinely independent, $S(C,R)$ be its touching sphere. Then $Q$ is nearly positively spanning if and only if $C \in \conv Q$. In this case $\defe Q = \sqrt{1-R^2}$.
\end{prop}

For further use we remark that in case $C\in \conv Q$ the touching sphere coincides with the circumsphere of $Q$.
By means of nearly positively spanning sets we may describe an  algorithm for F':

\begin{algo} \label{algo2}
\begin{enumerate}
\item
Let $j \in \{1,\dots,n\}$ be arbitrary, $x^1=a_j$, $Q_1=\{a_j\}$. 
\item
If $x^k/\|x^k\|$ is feasible, then stop with feasibility. Else let $m$ be the index of a most violated constraint for $x^k/\|x^k\|$, $y=x^k$.
\item
If $Q_k\cup \{a_m\}$ is positively spanning, then stop the algorithm.
\item
Compute the center $C$ of the touching sphere of $Q_k\cup\{a_m\}$.
\item
If $C\in  \conv (Q_k\cup\{a_m\})$, then let $x^{k+1}=C$, $Q_{k+1}=Q_k\cup\{a_m\}$, $k=k+1$. Goto step 2.
\item
Let $y$ be the point, where the line segment $\Ol{yC}$ intersects the relative boundary of $\conv (Q_k\cup\{a_m\})$. Let $F$ be a facet of $\conv (Q_k\cup\{a_m\})$, which contains $y$, and $a_j$ be the vertex of $Q_k$, which is not contained in $F$. Let $Q_k=Q_k\setminus \{a_j\}$ and goto step 4.
\end{enumerate}
\end{algo}

This algorithm shares a remarkable property with the simplex algorithm.

\begin{defi} \label{combi}
An algorithm for F, F' is said to be combinatorial, if the iterates depend only on subsets of the set of constraints.
\end{defi}

Thus for an instance of a combinatorial algorithm there are only finitely many possible iterates. As the iterates of Algorithm \ref{algo2} are centers of the touching sphere of $Q_k$, this algorithm is combinatorial like the simplex algorithm and different from {\sc Khachiyan}'s or {\sc Karmarkar}'s algorithm. To prove the convergence of a combinatorial algorithm it is clearly sufficient to show, that no iterate $x^k$ can be repeated.

\begin{lemma} \label{algo2_conv}
Algorithm \ref{algo2} solves F' for feasible as well as infeasible problems in finitely many steps.
\end{lemma}
Proof. By construction $\defe Q_k$ is strictly monotonely decreasing. \qed

It is certainly desirable to have some quantitative information on the progress of the algorithm. This can be obtained by a small change in Algorithm \ref{algo2}:

\begin{algo} \label{algo2m}
Construct the point $y$ in step 2 of the algorithm by choosing $y$
as the point on the line $\Ol{x^ka_m}$, which is closest to the origin.
\end{algo}

\begin{lemma}\label{progress}
Algorithm \ref{algo2m} solves the feasibility problem F' and it holds
\begin{equation} \label{prog_eq}
\defe Q_{k+1} \le \sqrt{\frac {1-v(x^k)^2}{1+(\defe Q_k)^2+2\,\defe Q_k\,v(x^k)}}\defe Q_k.
\end{equation}
\end{lemma}
Proof. The $y$ constructed in step 2 has the distance from the origin given by the right side of (\ref{prog_eq}) and by construction the distance of $x^{k+1}$ from the origin is certainly not greater than the distance of $y$ from the origin. \qed

While the estimate (\ref{prog_eq}) is not linear, it is independent of the dimension and in particular of the size of the problem. The formula may be used to estimate the number of steps, that the algorithm needs to reach a given level.

\begin{lemma}\label{speed}
In Algorithm \ref{algo2m} it holds $\defe Q_k \le t^{-1}$ after at most $t^2$ steps.
\end{lemma}
Proof. It is somewhat easier to consider instead of $\defe Q_k$ the quantities $y_k=(\defe Q_k)^{-1}$. We have from Lemma \ref{progress} and using $v(x^k)\ge 0$
\begin{equation} \label{Inv_prog}
y_{k+1}^2 \ge y_k^2+1.
\end{equation}
Together with $y_1=1$ this gives inductively $y^2_{k+1}\ge k+1$.
 \qed

It is not known (and was not investigated), whether the sequence of iterates in Algorithm \ref{algo2} and Algorithm \ref{algo2m} can differ for the same problem. According to Lemma \ref{progress} Algorithm \ref{algo2m} takes large steps, if the deficiency or the violation are large, in particular this is the case at the start of the algorithm.

We don't have information on the global behavior of the algorithm beyond Lemma \ref{speed}. We may obtain some more insight in its working by analyzing one step of Algorithm \ref{algo2}.  Using an appropriate system of coordinates we may assume that  $Q_k$ spans the affine plane $H=\{(x_1,0,\epsilon)\mid x_1\in \R^l\}$ for some fixed $\epsilon=\defe Q_k >0$ and appropriate $l$. Then we have $x^k=(0,0,\epsilon)$. Further we have that the constraint which is chosen in step 2 has the form $a_m=(x_1,x_2,\epsilon+\eta)$ with $x_1\in \R^l$, $x_2\in \R^{d-l}$ and $\eta < -\epsilon$. We find for the center $C$ of the touching sphere of  $Q_k\cup\{a_m\}$ in step 4 
$$
C=(0,0,\epsilon)+\frac{-\eta\epsilon}{\|x_2\|^2+\eta^2} (0,x_2,\eta).
$$

From this we find $\|C\| \approx \|x^k\|$, if $|\eta|\ll \|x_2\|$. In this case the algorithm can only increase the cardinality of $Q_k$, but not significantly decrease $\defe Q_k$. Clearly the crucial point of the algorithm is whether we may expect $C\in \conv (Q_k \cup \{a_m\})$. To see what happens, let us  denote by $a_m'$ the orthogonal projection of $a_m$ onto $H$. By construction we have that $x^k$ is the orthogonal projection of $C$ onto $H$. Finally let $y\in \relbd( \conv Q_k)$ be the point that $x^k \in \Ol {ya'_m}$. Now an examination of the triangle $y,a_m',a_m$ shows that $C\not \in \conv (Q_k \cup \{a_m\})$ if $\|y-x^k\|\ll \|a_m'-x^k\|$. Clearly this may happen and in a single step the algorithm may make very little progress. Fortunately this is not the complete story. 

In the algorithm we gather information about the $a_i$. If we define for $x\in S^d$
$$
M(x) = \{a\in S^d \mid a^T x\ge -v(x)\},
$$   
then we have for every iterate $x^j$ that $a_i \in M(x^j)$, $i=1,\dots,n$.
Consequently all $a_i$ satisfy for all $k$
$$
a_i \in \bigcap_{j=1}^k M(x^j).
$$
Unfortunately the sets $M(x^k)$ are not convex and this information is hard to use. 

If we want to repeat the bad step from above, the set $M(x^k)$ tells us that either for the new $\epsilon',\eta'$ holds $|\epsilon' +\eta'|\ll |\epsilon+\eta|$ or the angle between $a_m$ and $a_{m'}$ must be large. Thus we can hope that bad steps cannot be repeated too often. A precise analysis of the phenomena appears to be difficult.

\section{A polynomial algorithm}
In view of the last remark in the previous section it is near at hand to use periodically a rescaling to increase the deficiency in order to speed up the convergence. While it is simple to give a transformation of the sphere which takes polyhedra to polyhedra and increases the deficiency of a set $Q_k$ it appears difficult to work out the consequences of iterating such transformations.

Fortunately it turns out that for $\defe(Q_k)$ sufficiently small it is possible to find a transformation which increases the volume of the sets $P_A$. To work out our ideas it is necessary to recall some concepts of spherical convexity. A set $M \subset S^d$ is called convex, if there exists a half--sphere which contains $M$ and for every pair $x,y\in M$ there exists a great circle $C$ such that $C\cap M$ is connected. The intersection of convex sets is itself convex. By this we may define for every $M\subset S^d$ which is contained in a half--sphere the {\em spherical convex hull} $\co M$ as the intersection of all convex subsets of $S^d$ which contain $M$. Here we have chosen the notation $\co M$ to distinguish the spherical convex hull from the Euclidean convex hull used in the previous section. As every half--sphere is convex we have that the intersection of half--spheres is convex. We say that the intersection of finitely many half--spheres is a {\em (spherical) polyhedron}. An example of a spherical polyhedron is the set $P_A$ of feasible points of an instance $A$ of problem F'.

We denote for two points $x,y\in S^d$ the spherical distance by $\angle(x,y)$, i.e. $\cos (\angle(x,y))=x^Ty$. Next for $x\in S^d$ and $0\le \rho \le \pi$ 
$$
K_{\rho}(x) = \{y\in S^d \mid \angle(x,y)\le \rho\}
$$
denotes the {\em spherical cap} with center $x$ and radius $\rho$. For closed $M\subset S^d$ we may now define the spherical counterparts of inradius, circumradius and diameter. The {\em spherical inradius} $r'(M)$, the {\em spherical incenter} $c(M)$ respectively, are the radius, center respectively, of a largest spherical cap $K_{r'(M)}(c(M))\subset M$.  We remark that in contrast to the Euclidean case for convex $M$ the spherical incenter is always well defined. The {\em spherical circumradius} $R'(M)$, {\em spherical circumcenter} $C(M)$ respectively, are the radius, center respectively, of the smallest spherical cap $K_{R'(M)}(C(M))\supset M$. Here the notions $R',r'$ were chosen to distinguish spherical circumradius and inradius from their Euclidean counterparts. For $M\subset S^d$ we denote the boundary of $K_{R'(M)}(C(M))$ as {\em circumsphere} of $M$. For any $M$ contained in a half--sphere the {\em spherical diameter} $D(M)$ is given by
$$
D(M) = \sup \{\angle(x,y) \mid x,y\in M\}.
$$

For $M\subset S^d$, $a_1,a_2\in S^d$ satisfying $a_i^T x\ge 0$ for all $x\in M$, $a_i^T x=0$ for at least one $x\in M$ and $a_1\neq a_2$ we denote by $\pi-\angle(a_1,a_2)$ a {\em width} of $M$. The {\em breadth} $d(M)$ is then defined as the minimum over all widths of $M$.

For $M\subset S^{d}$ the {\em polar set} $M^*$ is defined by
$$
M^* = \{x \in S^d \mid m^Tx\ge 0 \text{ for all $m\in M$}\}.
$$

In our context it is of importance to observe that the polar set of a spherical polyhedron is again a polyhedron and that inclusions are reversed, i.e. for $M \subset M'$ we have $M^* \supset {M'}^*$. For $M=\{x^1,\dots,x^k\}$, $x^1,\dots,x^k$ linearly independent, we have that $M^*$ is a polyhedron with $k$ facets, which we denote as a $k-1$--cosimplex.

Further for $M\subset S^d$ convex we have the relations
$$
r'(M)+R'(M^*) = \pi, \qquad d(M)+D(M^*)=2\pi.
$$
    
In this terminology we have for an instance $A$ of F' and set $Q=\{a_{j_1},\dots,\linebreak[2] a_{j_l}\}$ from step 5 of Algorithm \ref{algo2} with $P_A\neq \emptyset$ that $M=\co Q$ is an $l-1$--simplex with $M\subset P_A^*$ and $\cos(R'(M)) =\defe Q$. Correspondingly $M^*$ is a $l-1$--cosimplex with $P_A \subset M^*$ and $r'(P_A)\le r'(M^*)=\pi-R'(M)$.

Moreover for  $M\subset S^d$ the quantities $R'(M)$ and $D(M)$ are not independent as the following theorem by {\sc L. A. Santal\'o} \cite {S} shows:

\begin{satz}[\sc Santal\'o]
Let $M\subset S^d$ be contained in a half--sphere. Then
\begin{enumerate}
\item
For $\cos R' \ge 1/\sqrt {d+1}$ 
$$
\cos 2R' \le \cos D \le \frac{(d+1)\, \cos^2R'-1}{d}.
$$
\item
For $0\le \cos R' \le 1/ \sqrt{d+1}$
$$
\cos 2R' \le \cos D \le \frac{(d+1)\cos^2 R'-1}{1+(d+1)\cos^2 R'}\quad \text{for $d$ odd},
$$
\begin{multline*} 
\cos 2R' \le \cos D \le
\frac{(d+1)\,\cos^2 R' -1} {(1+(d+1)\cos^2 R')^{1/2}\left(1 +\frac{(d+1)(d-2)}{d+2}\cos^2 R'\right)^{1/2}}\\ \text{for $d$ even}.
\end{multline*}
\end{enumerate}
\end{satz}

Using {\sc Santal\'o}'s theorem we may estimate the breadth $d(P_A)$ of $P_A$ by computing $R'(Q)$. After a suitable transformation of coordinates we have for all $(x_1,\dots,x_{d+1})^T\in P_A$ that $|x_1|\le \sin d(P_A)$, i.e. all points of $P_A$ are close to the equator $\{(x_1,\dots,x_{d+1})\in S^d \mid x_1=0\}$, and $P_A$ can be enlarged by a transformation. If we know that in the beginning  $P_A$ is not to small, we must find a feasible $x^k$ after a number of steps, for otherwise the transformed $P_A$ would not fit on the sphere.

A suitable quantity to measure the size of $P_A$ is the volume and the estimate of the size of $P_A$ at the start of the algorithm is standard from the theory of linear programming. This leaves the determination of the hyperplane. While the proof of {\sc Santal\'o}'s theorem is constructive, the construction of the hyperplane needs the computation of about $\binom{d}{d/2}$ distances and is therefore not practicable for a polynomial algorithm. Thus we substitute {\sc Santal\'o}'s theorem by a weaker one, for which we can easily compute all quantities.

\begin{defi}\label{Eck_durch}
Let $a_1,\dots,a_{d+1}\in S^d$, such that 
 $Q=\co\{a_1,\dots,a_{d+1}\}$ 
 is a spherical simplex. For each $a_i$ let $a_i'$ be the intersection of the circle through $a_i$ and the circumcenter $C(Q)$ with $\co\{a_1,\dots,a_{i-1},a_{i+1},\dots,a_{d+1}\}$. Then the {\em vertex diameter} $D_v(Q)$ is defined by
$$
D_v(Q) = \max\{\angle(a_1,a'_1),\dots,\angle(a_{d+1},a'_{d+1})\}.
$$
\end{defi}

Clearly the vertex diameter can be computed easily. While \cite{S} contains no explicit estimate of the vertex diameter, the methods can be adopted. To be complete we present the changes in the proof.
First we have for the vertex diameter of the regular simplex
\begin{lemma}[Santal\'o] \label{Reg_Eck_durch}
Let $Q\subset S^d$ be a regular simplex with circumradius $R'$ and vertex diameter $D_v$, then
\begin{equation} \label{Reg_ver}
\cos D_v = \frac{(d+1)\,\cos^2 R'-1}{(1+(d-1)(d+1)\cos^2 R')^{1/2}}.
\end{equation}
\end{lemma}
Proof. The formula is a simple combination of formulas (2.3) and (2.13) in \cite{S}. \qed

\begin{lemma}[Santal\'o]\label{circum_for}    
Let $Q=\co\{a_1,\dots,a_{d+1}\}\subset S^d$ be a spherical simplex, such that all vertices are contained in its circumsphere. Let $C=\sum_{i=1}^{d+1} \mu_i a_i$ be the circumcenter of $Q$. Then for the circumradius $R'$ of $Q$ it holds
$$
\cos R' = \left(\sum_{i=1}^{d+1} \mu_i\right)^{-1},\qquad
\frac{1}{\cos^2 R'} \le (d+1)\sum_{i=1}^{d+1} \mu_i^2.
$$
\end{lemma}
Proof. We have
$$
\sum_{i=1}^{d+1} \mu_i\, \cos R' = \sum_{i=1}^{d+1} \mu_i\,a_i^T C = C^TC=1.
$$
From the equation we obtain
$$
\frac{1}{\cos^2 R'} = \sum_{i=1}^{d+1} \mu_i^2 +2 \sum_{1\le i<j\le d+1}\mu_i\mu_j.
$$
The inequality is an immediate consequence of $2\mu_i\mu_j \le \mu_i^2+\mu_j^2$. \qed    

\begin{lemma}[Santal\'o]\label{Eck_durch_Ab}
Among all spherical $d$--simplices on $S^d$ with circumradius $R'$ and $\cos R' <1/\sqrt{d+1}$, which contain all vertices in their circumsphere, the regular simplex has minimal vertex diameter.
\end{lemma}
Proof. Assume that there is a simplex $Q\subset S^d$ with circumradius $R'$, circumcenter $C$ and vertex diameter $\Delta_v$, which has all vertices on its circumsphere and has smaller vertex diameter than the regular simplex with circumradius $R'$. Let $D_v$ be the vertex diameter of the regular simplex with circumradius $R'$. From (\ref {Reg_ver}) we have $\cos D_v<0$. Let $a_1,\dots,a_{d+1}$ be the vertices of $Q$. We have
$$
C= \sum_{i=1}^{d+1} \mu_i a_i \quad \text{with} \quad \mu_i>0.
$$
and
$$
a'_1 = \sum_{i=2}^{d+1}\mu_i a_i/\|\sum_{i=2}^{d+1}\mu_i a_i\|.
$$
From this we obtain
\begin{equation} \label{Mittel}
\begin{split}
\mu_1 \cos R' &= \mu_1 a_1^T C  \\
&= \mu_1^2 + \mu_1 \left(\sum_{i=2}^{d+1} \mu_i^2+2\sum_{2\le i<j\le d+1} \mu_i\mu_j a_i^T a_j\right)^{1/2} \cos \angle(a_1,a_1').
\end{split}
\end{equation} 

By our assumption we have $ \cos \angle(a_1,a_1')>\cos D_v$. Thus (\ref{Mittel}) yields  
\begin{equation}\label{Weiter}
\mu_1 \cos R'  - \mu_1^2 - \mu_1 \left(\sum_{i=2}^{d+1} \mu_i^2+2\sum_{2\le i<j\le d+1} \mu_i\mu_j a_i^T a_j\right)^{1/2} \cos D_v > 0.
\end{equation}

For $a_2,\dots,a_{d+1}$ we obtain inequalities analogous to (\ref{Weiter}). Summing up all these inequalities, using the Cauchy--Schwarz inequality $x^Ty \le \|x\|^{1/2}\|y\|^{1/2}$ and taking $\|C\|=1$ into account we obtain
\begin{equation} \label{Lang}
\begin{split}
0 &<
\sum_{k=1}^{d+1} \mu_k \cos R' - \sum_{k=1}^{d+1} \mu_k^2 
-\sum_{k=1}^{d+1} \mu_k\left( \sum_{\substack{i=1\\ i\neq k}}^{d+1} \mu_i^2+2\sum_{\substack{1\le i<j\le d+1\\ i,j\neq k}} \mu_i\mu_j a_i^T a_j\right)^{1/2}\cos D_v\\
&<\sum_{k=1}^{d+1} \mu_k \cos R' - \sum_{k=1}^{d+1} \mu_k^2-\\
&\quad \left( \sum_{k=1}^{d+1} \mu_k^2\right)^{1/2}\left(d \sum_{i=1}^{d+1} \mu_i^2 +(d-1)\left(1-\sum_{i=1}^{d+1} \mu_i^2\right)\right)^{1/2}\cos D_v\\
&=\sum_{k=1}^{d+1} \mu_k \cos R' - \sum_{k=1}^{d+1} \mu_k^2-
\left( \sum_{k=1}^{d+1} \mu_k^2\right)^{1/2}\left( \sum_{i=1}^{d+1} \mu_i^2 +(d-1)\right)^{1/2}\cos D_v
\end{split}
\end{equation}

Solving (\ref{Lang}) for $\cos D_v$ and using the identity and inequality from Lem\-ma~\ref{circum_for} gives finally
\begin{align*}
\cos D_v &< \frac{1-\sum_{k=1}^{d+1} \mu_k^2}{\left( \sum_{k=1}^{d+1} \mu_k^2\right)^{1/2}\left( \sum_{i=1}^{d+1} \mu_i^2 +(d-1)\right)^{1/2}}\\
&\le \frac{(d+1) \cos^2 R' -1}{\left(1+(d+1)(d-1)\cos^2 R'\right)^{1/2}}
\end{align*}

This is a contradiction to (\ref{Reg_ver}). \qed  
    
By Lemmas \ref{Reg_Eck_durch}, \ref{Eck_durch_Ab} we may determine with the help of a set $Q$ two half--spheres, which contain $P_A$ in their intersection and the angle $\phi$ between them is most
\begin{equation}\label{min_phi}
\cos \phi \ge \frac{1-(d+1)\cos^2 R'(Q)}{(1+(d-1)(d+1)\cos^2 R'(Q))^{1/2}}.
\end{equation}

(\ref{min_phi}) yields by an easy calculation and further estimation 
\begin{equation}\label{min_sin_phi}
\sin \phi \le (d+1) \cos R'.
\end{equation}

Next we show that for $\phi$ sufficiently small we may find a transformation of the sphere, which maps spherical polyhedra into spherical polyhedra and increases the volume of $P_A$ by a prescribed factor.
Any spherical map
$$
B: S^d \to S^d, \quad B(x)=B'x/\|B'x\|
$$
with a nondegenerate linear map  $B'$ has the first property. The influence of $B$ on the volume is most easily studied by use of parametrisations.  

A parameterized $d$--surface in $\R^{d+1}$ is a differentiable injective map $\Phi:G\subset \R^{d} \to \R^{d+1}$ such that for every $x\in G$ the Jacobi--matrix $\left( \frac {\partial \Phi}{\partial x}\right)$ is nondegenerate. For $M\subset \Phi(G)$ the volume of $M$ is given by
$$
V(M) = \int_{\Phi^{-1}(M)} \det \Phi\, dx,
$$
where $\det \Phi$ denotes the volume of the $d$--dimensional parallelepiped spanned by $\frac{\partial \Phi}{\partial x_1},\dots,\frac{\partial \Phi}{\partial x_d}$. For a differentiable map $\psi : \R^{d+1}\to \R^{d+1}$ with appropriate properties $\Psi \circ \Phi$ is a parameterized surface and consequently
$$
V(\Psi(M)) = \int_{\Phi^{-1}(M)} \det(\Psi\circ \Phi)\,dx
$$
is the volume of $\Psi(M)$.
Thus $(\det(\Psi\circ \Phi)/\det \Phi)(\Phi^{-1}(x))$ is the local change of volume in $x\in \Phi(G)$ under the map $\Psi:\Phi(G)\to \R^{d+1}$. 

Here we are interested in the case
\begin{gather}
\Phi: B^d \to \R^{d+1}, \quad \Phi(x_1,\dots,x_d)=(x_1,\dots,x_d,\sqrt{1-\sum_{i=1}^d x_i^2})^T,\\
\begin{split}\label{Psi}
\Psi_\alpha&:\R^{d+1} \to \R^{d+1},\\ \Psi_{\alpha}&(x_1,\dots,x_{d+1})= (\alpha x_1,x_2,\dots,x_{d+1})^T/\sqrt{1+(\alpha^2-1)x_1^2}, \quad \alpha>0.
\end{split}
\end{gather}

The restriction of $\Psi_\alpha$ on $S^d$ is the spherical transformation $(\alpha x_1,x_2,\dots,\linebreak[2]x_{d+1})^T/\|(\alpha x_1,x_2,\dots,x_{d+1})^T\|$. We have chosen $\Psi_\alpha$ rather than the other transformation, as its  derivatives can  be computed more easily.  In particular we consider the point $\Ol x=(\beta,0,\dots,0,\sqrt{1-\beta^2})^T$ with $0\le \beta < 1$. Computation of the Jacobi matrices and of the determinants with the help of Gram determinants yields
$$
\det \Phi(\Ol x) = \frac 1 {(1-\beta^2)^{1/2}}, \quad \det (\Psi _\alpha \circ \Phi)(\Ol x)=\frac{\alpha}{(1-\beta^2)^{1/2} (1+\beta^2(\alpha^2-1))^{(d+1)/2}}.
$$
Thus the transformation $\psi_\alpha$ gives for all $(x_1,\dots,x_{d+1})^T\in S^d$ with $x_1=\beta$ the local change of volume 
$$
f(\beta)=\frac{\alpha}{(1+\beta^2(\alpha^2-1))^{(d+1)/2}}.
$$
Consequently we obtain for $M\subset S^d$ with $x_1\le \beta$ for all $(x_1,\dots,x_{d+1})^T\in M$ and $\alpha>0$ the estimate
\begin{equation} \label{Vol_Wechsel}
V(\Psi_\alpha(M)) \ge\frac{\alpha}{(1+\beta^2(\alpha^2-1))^{(d+1)/2}}V(M).
\end{equation}

We may summarize the previous considerations in
\begin{lemma} \label{Vol_Trans}
Let $M\subset S^d$ such that for all $(x_1,\dots,x_{d+1})^T\in M$ it holds $x_1\le \beta$. Further let $\Psi_\alpha : S^d\to S^d$, $\alpha>0$, be given by (\ref{Psi}). Then the volume of $\Psi_\alpha(M)$ can be estimated by (\ref{Vol_Wechsel}).
\end{lemma}
  
Evidently we may generalize the transformation $\Psi_\alpha$ to the case of a transformation orthogonal to an arbitrary equator $a^Tx=0$ instead of the equator $x_1=0$.

As we work with the homogenized problem, we have to study the influence of (\ref{Homogen}) on the volume. Using the same technique as before we find for $M \subset \R^d$ for the volume of the homogenized set
$$
V(\Phi(M)) = \int_M \frac{1}{(\|x\|^2+1)^{(d+1)/2}}dx.
$$
Thus we have
\begin{lemma}\label{Homo_vol}
Let $M\subset \R^d$ with $\|x\|\le R$ for all $x\in M$. Then $V(\Phi(M))\ge (R+1)^{-(d+1)}$.
\end{lemma}

The following considerations concerning polynomial algorithms could be done in the context of oracle--polynomial time algorithms for convex bodies in \cite{GLS}. The presentation becomes somewhat simpler if we follow the presentation in \cite{PS} for algorithms for the feasibility problem of linear systems.

According to the corollary to Lemma 8.7 in \cite{PS} we can solve F in polynomial time, if we can solve the following problem in polynomial time: For the set
\begin{equation}\label{Ungleich}
M= \{x\in \R^d \mid a_i^Tx >b_i,\; i=1,\dots,n\}
\end{equation}
of size $L$ find a point in $M$ or show that there exists no such point in  time polynomial in $L$. Here the size $L$ of the system measures the amount of space needed to write down the inequalities. For a precise definition of the size compare \cite{GLS} or \cite{PS}.  

We have that either  $M$ is empty or the volume of $M$ is not too small:
\begin{lemma}\label{Min_Gro}
If the set $M$ from (\ref{Ungleich}) is nonempty, then
$$
V(M\cap \{x\mid \|x\|\le d\,2^L\}) \ge 2^{-(d+2)L}.
$$
\end{lemma}
Proof. \cite{PS}, Lemma 8.14. \qed

For the $a_i,b_i$ from (\ref{Ungleich}) let $a'_i=(a_i^T,-b_i)^T$. Then we have for the homogenization of $M$
$$
\Phi(M) = \{x\in S^d \mid {a'_i}^Tx>0,\;i=1,\dots,n,\; x_{d+1}>0\}.
$$

We may combine Lemmas \ref{Homo_vol}, \ref{Min_Gro} to obtain
\begin{lemma} \label{Homo_min}
Let $M$ be the set from (\ref{Ungleich}) with size $L$. Then
$$
V(\Phi(M)) \ge 2^{-3(d+2)L}
$$
or $\Phi(M)$ is empty.
\end{lemma} 

Now we have collected the parts to show that the following modification of Algorithm \ref{algo2m} is polynomial.

\begin{algo} \label{poly_alg}
Let $\beta_d = \sqrt{((4/3)^{2/(d+2)}-1)/3}$, $max = 6(d+2)L$.

Main loop:
\begin{enumerate}
\item
Let $p=0$.
\item
Let $p=p+1$. If $p>max$, stop with infeasibility.
Let $j\in \{1,\dots,n\}$ be arbitrary, $k=1$, $x^1=a_j$, $Q_1=a_j$.
\item
If $x^k/\|x^k\|$ is feasible, then stop with feasibility. Else let $m$ be the index of any constraint, which is violated for $x^k/\|x^k\|$.  If $Q_k \cup\{a_m\}$ is a positively spanning set, then stop with infeasibility. Let $y$ be the point on the line segment $\Ol{x^ka_m}$, which is closest to the origin.
\item
Compute the center $C$ of the touching sphere of $Q_k\cup\{a_m\}$.
\item
If $C \in  \conv (Q_k\cup \{a_m\})$, then do the following: Let $x^{k+1}=C$, $Q_{k+1}=Q_k\cup\{a_m\}$, $k=k+1$. If $\defe(Q_{k+1})<\beta_d/(d+1)$ then Transform and goto step 2 else goto step 3.
\item
Let $y$ be the intersection of the line segment $\Ol{yC}$ with the relative boundary of $\conv (Q_k \cup\{a_m\})$. Let $F$ be a facet of $\conv (Q_k \cup \{a_m\})$, which contains $y$. Let $a_j$ be the vertex of $\conv (Q_k \cup \{a_m\})$, which is not contained in $F$. Let $Q_k=Q_k\setminus\{a_j\}$. Goto step 4.
\end{enumerate}

Transform:
\begin{enumerate}
\item
Determine the vertex diameter of $\co Q_k$. Let $a_i$ be the vertex which gives the vertex diameter.
\item
Apply the transformation $\Psi_2$ with respect to the equator $a_i^Tx=0$.
\end{enumerate}
\end{algo}

\begin{satz}\label{Haupt}
Algorithm \ref{poly_alg} solves F in polynomial time.
\end{satz}
Proof. We denote by $M_p$ the set of feasible points before the $p$--th transformation. First we show the correctness of the algorithm. The first stopping condition in step 3 gives trivially the correct result. If the algorithm stops at the second stopping condition, then $M\subset P_{Q_k}$. $P_{Q_k}$ is a lower dimensional sphere and thus $V(P_{Q_k})=0$ and the instance is infeasible by Lemma \ref{Homo_min}.

Now we consider the state of the algorithm at the call of the procedure ``Transform''. We have $\cos R'(Q_k) \le \beta_d/(d+1)$. By (\ref{min_sin_phi}) we find that all points of $M_p$ are contained in the zone
$\{x\in S^d \mid 0\le a^T_i x \le \beta_d\}$ for the $a_i$ giving the vertex diameter. From this we obtain by Lemma \ref{Vol_Trans} that $V(M_{p+1})/V(M_p)\ge 3/2$. Finally we obtain by Lemma \ref{Homo_min} that after $\max$ steps we have $V(M_{max})\ge 2 > V(S^d)$. Consequently an existing feasible point must be found before this step. This proves the correctness of the algorithm.

To calculate the number of computations we only have to estimate the number of iterations in the inner loop starting in step 2. By a power series expansion we see easily that $\beta\ge c\,d^{-1/2}$ for suitable $c>0$. Thus by Lemma \ref{speed} there are $O(d^3)$ iterations.

It remains to show that it is sufficient to do all calculations with rational numbers whose size is bounded by a polynomial in $L$. This can be done in the same way as in {\sc Khachiyan}'s algorithm and we skip the details. \qed 

Some remarks concerning the algorithm:
\begin{enumerate}
\item
If we compare Algorithm \ref{algo2m} and Algorithm \ref{poly_alg} we find that we have replaced the condition  in Algorithm \ref{algo2m} that the new $a_m$ in step 2 is given by  a most violated constraint by the weaker condition that the new $a_m$ in step 4 of Algorithm \ref{poly_alg} is given by any violated constraint. We can easily check that it is sufficient to choose any $a\neq 0$, such that $a^T x^k<0$, but $a^T x\ge 0$ for all $x\in P_A$, i.e. if $\{a^T x=0\}$ is a separating great sphere for $x^k$ and $P_A$. Moreover we made no use of the fact that $P_A$ is a polyhedron, but used only that the volume of $P_A$ cannot be too small. Now a separating hyperplane for the unhomogenized problem immediately yields a separating great sphere for the homogenized problem. Thus the algorithm is essentially a polynomial algorithm for the following convex feasibility problem: Given is a convex set $K\in \R^d$ which is either empty or well bounded in the language of \cite{GLS} and an oracle which for any $x\in \R^d$ confirms that $x\in K$ or gives a separating hyperplane. Under these assumptions find a point in $K$  or show that $K$ is empty.
\item
While the previous remark pointed out that from a theoretical point of view, it is sufficient to choose any separating great sphere in step 4, from a practical point of view this is certainly not the case. This is shown by the following example, which shows that the inner loop in Algorithm \ref{poly_alg} may need many iterations:

In this example the sets $Q_k$ have always 3 elements which are for simplicity unnormalized. Let
$$
Q_k = \{(1,\epsilon,0,\dots,0)^T,(-1,\epsilon,0,\dots,0)^T,(0,-\delta_k,1,0,\dots,0)^T\},
$$
where $\epsilon$ and $\delta_k$ are positive numbers of the kind that $\delta_k$ is small with respect to 1, and $\epsilon$ is small with respect to $\delta_k$. The points $x^k$ are of the form $x^k=(0,1,\eta_k,0,\dots,0)^T$. $\{(0,-\delta_{k+1},1,0,\dots,0)^Tx=0\}$ is a separating great sphere if $\delta_{k+1}>\eta_k$. The computation of $x^{k+1}$ gives
\begin{equation}\label{neu_x}
x^{k+1}=(0,1,\eta_{k+1},0,\dots,0)^T, \quad
\eta_{k+1}= \delta_{k+1}+\epsilon \frac{\sqrt{1+\delta_{k+1}^2}}{\sqrt{1+\epsilon^2}}.
\end{equation}
Thus if $\delta_{k+1}$ is close to $\eta_k$ we have that the distance $\angle(x^k,x^{k+1})$ is of order $\epsilon$. Now the iteration  may be repeated with a $\delta_{k+2}$ which is of the same order as $\delta_{k+1}$. This behavior does not affect the polynomiality as $\epsilon$ is bounded below by a function depending on the dimension, but demonstrates that the convergence may be rather slow.

\item
If we choose in step 4 the most violated constraint then this phenomenon cannot occur, as is shown by the discussion after Lemma \ref{speed}. 


\item
While the choice of $x^k$ as circumcenter of the sets $Q_k$ is fundamental for our work, this is certainly not the case for the construction of the transformation. Although it did secure the polynomiality of the algorithm, there are some apparent disadvantages. We have to work with three different quantities (circumradius, vertex diameter, volume) and these quantities are rather weakly related. Moreover the volume does not behave very nicely for transformations of the sphere. The apparent choice for a transformation would be to transform after each iteration $\co Q_k$ into a regular simplex of the same circumradius such that the circumcenter of $\co Q_k$ is transformed to the circumcenter of the resulting regular simplex. Unfortunately we don't know of any quantity which behaves nicely under this transformation and thus don't even know, whether the resulting algorithm is convergent.

\item
Still there are probably other transformations which are easier to handle than the one sketched above and lead to polynomial algorithms. One approach might be the following: We want to construct in polynomial time a positively spanning set for infeasible problems. For our heuristic considerations we assume that this positively spanning set is actually a positive basis. As we want to construct a positive basis, we look at the violation. Considering the fact that for a regular simplex the quotient of the circumradius and the inradius is $d$ we find that $v(x^k)\ge 1/(d+1)$, if the points of the positive basis form a regular simplex.  Now, let us assume that we have found a point $x^k$ such that $v(x^k)$ is small compared to $1/(d+1)$. Let $Q$ be an arbitrary positive basis spanned by the $a_i$. Then it is easy to see, that there must be a face $F$ of $\conv Q$, such that $F$ is nearly parallel to $\{(x^k)^T x=0\}$ and has a small distance from the origin. Let us  apply for suitable $\alpha>1$ a transformation $\Psi_\alpha$ with respect to the equator $\{(x^k)^Tx=0\}$. Locally we increase the violation. Globally it may be seen, that the distance to the origin of all faces of $\conv Q$, which have a dimension less or equal than the dimension of $F$ is increased. One may hope that the minimal rate of increase is given by a function, which depends only on the dimension. Thus one should obtain fast, i.e. polynomial convergence.

We have not tried to work out the theoretical details of this algorithm. Thus we don't know whether it works and have no information on its theoretical quality. We may say that this algorithm would avoid the problems of our algorithm which we pointed out in the previous remark, as in our context it is certainly more natural to deal with distances than with volume. Moreover the transformation with a pole $x^k$ is more natural than the transformation with a pole $a_i$ and  it is easy to implement this algorithm. This was in fact done with a transformation taking place whenever $v(x^k)\le1/\sqrt d$. The surprisingly good results of this algorithm are presented in the next section.  
\end{enumerate}

Altogether we may say, that solely as a consequence of the weak estimate in Lemma \ref{speed} it is possible to construct a sequence of transformations which transform the combinatorial Algorithm \ref{algo2m} into a polynomial one. This is a remarkable contrast to the combinatorial simplex algorithm, as for the simplex algorithm it is even unknown, whether an optimal choice of pivot--steps would lead to a polynomial algorithm ({\sc Hirsch}--conjecture, for an account see \cite{Z}). 

Moreover there is some hope, that Algorithm  \ref{algo2} is polynomial. That would  be of particular interest, as an algorithm which is combinatorial as well as polynomial might be a candidate for a strongly polynomial algorithm for F. It is still an open question whether such an algorithm exists. For  details about this problem see \cite{GLS}.

\section {Implementation of the algorithms and numerical experiments}
Up to now we have only studied the geometry of the algorithms and we still have the problem to perform the necessary calculations. We don't want to go into details or find the optimal way to do this but just show that there are no particular problems and that the complexity of a single step is on average comparable to the simplex algorithm, i.e. $O(dn)$.

The only steps in Algorithm \ref{algo2} where the computations are not completely obvious are step 3 and step 4. Here we must compute the center $C$ of the touching sphere and check whether it is in the convex hull of the $a_i$. We observe that by our construction  $Q_k$ is linearly independent and $Q_k\cup\{a_m\}$ affinely independent. We compute a least squares representation
\begin{equation}\label{least_squares}
a_m=a+\sum_{a_i\in Q_k} \lambda_i a_i, \quad a \text{ orthogonal to all }a_i\in Q_k
\end{equation}
and distinguish two cases.

\begin{enumerate}
\item 
$a=0$. Then $Q_k\cup\{a_m\}$ is linearly dependent, $C=0$ and using (\ref{least_squares}) we easily find a representation of $0$ as an affine combination of $Q_k\cup\{a_m\}$.
\item
$a \neq 0$. We observe that $C$ is characterized by the following three properties:
\begin{enumerate}
\item  $a_i^T C=a_j^T C>0$ for all $a_i,a_j\in Q_k\cup\{a_m\}$. 
\item  $C=\sum_{a_i\in Q_k\cup\{a_m\}} \mu_i a_i$ with $\mu_i \in \R$.
\item $\sum \mu_i=1$.
\end{enumerate}

Thus if a $C'\neq 0$ has the first two properties, then $C$ is just a multiple of $C'$.

 Let us look at the points $y$ constructed in step 6 of the algorithm.  We find inductively that we have
\begin{gather*}
\angle(y,a_i)=\angle(y,a_l)\le \angle(y,a_m) \quad \text{for all $a_i,a_l\in Q_k$}\\
y = \sum_{a_i\in Q_k\cup\{a_m\}} \mu_i a_i, \quad \text{with all $\mu_i\ge 0$}.
\end{gather*}

Thus using (\ref{least_squares}) we find that there is a $C'$ of the form $C=y+\rho a$ for suitable $\rho>0$. From this we determine $C$ as an affine combination of $Q_k$.
\end{enumerate}

The next $a_i$ to leave $Q_k$ is then apparently found by a simple line search. It remains to solve (\ref{least_squares}). To do this we consider the $Q_k$ as matrices with columns $a_i$. We observe that consecutive matrices $Q_k$, $Q_{k+1}$ only differ in that a column is added or taken away. If we solve (\ref{least_squares}) by computing a QR--factorization of $Q_k$ then it is well known that the QR--factorizations of consecutive $Q_k$, $Q_{k+1}$ can be obtained by an updating process which takes only $O(d^2)$ computations, cf.~\cite{GL}. Thus taking into account that the determination of $a_m$ in step 2 takes $O(dn)$ computations and we have on average to add one column to $Q_k$ and take one column away, we find that a step of Algorithm \ref{algo2} takes $O(dn)$ computations on average.

Clearly our results above give no information whether the algorithms could be of practical use. An obvious way to test this would be to solve the test problems as provided by http://www.netlib.org/lp. These problems are sparse optimization problems, i.e. to solve them we have to deal with many issues which are rather unrelated to the basic algorithm, but will strongly influence the results. Among these are solving the optimization problem by an algorithm for feasibility, dealing with equality--constraints, dealing with the sparsity. Thus we have taken a much more simple approach. We generated three different sets of problems. The most basic way was done by generating in a random way $n$ vectors of unit length $a_i$ and $n$ negative constants $b_i$ (Ex 1 in the following tables). Thus we obtain a system of inequalities $a_i^Tx\ge b_i$ which has the origin as a feasible point. In the next problem we  chose $b_i=0$ for $i=1,\dots,d+1$ and $a_{d+1}=-\sum_{i=1}^d a_i/\|\sum_{i=1}^d a_i\|$. Thus in this case we obtained with probability 1  problems which have exactly one feasible point (Ex 2). The last set differed from the second only in that we chose $b_{d+1}=10$ to construct infeasible problems (Ex 3). In order to use the origin as a starting point all systems were translated by a random vector.

If we look at the problems then we see that after choosing the correct set of indices the second problem is equivalent to the solution of a system of linear equations. If we use the Gauss--algorithm, and allow to transform the complete problem in every step then we see that a solution would involve $d$ steps, each with $O(nd)$ operations. In particular the number of steps would be independent of the number of inequalities. Certainly we would consider such an algorithm as excellent from a practical point of view and we may ask how close we can get.

We used three different algorithms for the solution of the problems. First we used the following variant of the simplex algorithm to have a basis of comparison. We consider the problem as the dual of a problem in standard form. We perform $d$ pivot--steps to generate a basis. This corresponds to ask for equality in the inequalities corresponding to the basic columns. Then we choose a positive constant right side and use the simplex--algorithm to solve the generated problem. If the algorithm stops with a finite solution, this shows the existence of a feasible point in the dual problem in the intersection of the constraints corresponding to the basic columns. If the algorithm stops with unboundedness this shows the infeasibility of the dual problem. This algorithm is denoted Alg 1 in the following  in the tables.

The next algorithm was a straightforward implementation of Algorithm \ref{algo2} (Alg 2). The last experimental algorithm (Alg 3) used an additional rescaling in Algorithm \ref{algo2}. Here it does not make much sense to use the polynomial Algorithm \ref{poly_alg} as the number of steps before the first transformation in step 5 takes place is not much smaller than the total number of steps for the complete solution by Algorithm \ref{algo2}. Instead we used the following procedure. We have by (\ref{prog_eq}) that every step of the algorithm makes relative progress of at least $\sqrt{1-v(x^k)^2}$. Thus when $v(x^k)$ becomes less than  $1/\sqrt d$ then we do  the following. Let $v(x^k)$ be determined by the constraint $a_r$, i.e. $v(x^k)=-a_r^T x^k/\|x^k\|$. We set $a_i'=(I+\lambda x^k {x^k}^T)a_i$, $i=1,\dots,n$ and renormalize the $a_i'$ to unit length. Here $\lambda$ is chosen in a way, such that $a_r' x^k/\|x^k\|=-\sqrt{2/ d}$. We observe that the matrix $Q_k$ is transformed by the formula   $Q_k=\mu (Q_k+w\, {x^k}^T)$ for a scalar $\mu$ and a vector $w$ which are easily determined. For such a transformation it is possible to update the QR--factorization  in $O(d^2)$ computations, cf \cite{GL}. Thus the transformation takes a total of $O(nd)$ computations.

To test the dependence on the dimension we chose the dimensions $d=10$, $20$, $40$, $80$, $160$, $320$ and $640$ and $8d$ inequalities for each dimension. To test the dependence on the  number of inequalities we chose $d=100$ and $n=400,800,1600,3200,6400$. For each case 5 examples were generated and solved. The average numbers of steps are shown in the tables.
\bigskip

 


\begin{table}[htb] {\centering{\small
\begin{tabular}{|c|r|rrrr|}\hline
&\multicolumn{1}{|c}{dim} &\multicolumn{1}{|c} {Alg 1} & \multicolumn{1}{c}{Alg 2} & \multicolumn{2}{c|}{Alg 3}\\ \hline
&10&13.4 & 20.6 & 16.8 & 5.8\\
&20&31.2 & 57.8 & 35.2 & 12.0\\
&40& 93.4 & 146.0 & 69.2 & 18.8\\
Ex1&80 & 274.4 & 353.2 & 146.6 & 23.2\\
&160 & 769.2 & 926.0 & 294.8 & 27.4\\
&320 & 2114.6 & 2156.6 & 585.0 & 29.8\\
&640 & 6321.0 & 4756.4 & 1179.0 & 33.8\\ \hline
&10 &17.0 & 28.8 & 25.8 & 9.2\\
&20& 42.6 & 62.6 & 54.2 & 14.0\\
&40 & 122.0 & 154.8 & 108.8 & 20.8\\
Ex2&80 & 346.6 & 385.2 & 228.8 & 29.6 \\
&160& 925.0 & 923.8 & 528.2 & 35.0\\
&320& 2528.0 & 2296.8 & 939.0 & 41.2 \\
&640 & 7294.4 & 5388.0 & 1909.4 & 45.4\\ \hline
&10& 15.8 & 27.4 & 24.0 & 7.4 \\
&20 & 41.0 & 58.4 & 50.2 & 11.0 \\
&40 & 112.4 & 141.2 & 101.2 & 17.0 \\
Ex3&80& 326.2 & 368.6 & 210.0 & 19.8 \\
&160& 887.0 & 857.2 & 422.0 & 21.0 \\
&320 & 2565.2 & 2183.2 & 867.0 & 22.6 \\ 
&640 & 7151.0 & 5125.2 & 1787.0 & 23.4\\
\hline
\end{tabular}}
 
\caption{
Number of steps in dependence of the dimension}

}   
\end{table}

\begin{table}[htb]
\centering {\small
\begin{tabular}{|c|rr|rr|rr|rr|} \hline
&\multicolumn{4}{|c|}{Alg 1} & \multicolumn{2}{c|}{Alg 2}&\multicolumn{2}{c|}{Alg 3}\\
&\multicolumn{1}{|c}{$\alpha$}&\multicolumn{1}{c|}{$\beta$}&
\multicolumn{1}{|c}{$\alpha$}&\multicolumn{1}{c|}{$\beta$}&
\multicolumn{1}{|c}{$\alpha$}&\multicolumn{1}{c|}{$\beta$}& 
\multicolumn{1}{|c}{$\alpha$}&\multicolumn{1}{c|}{$\beta$}\\ \hline
Ex1 &0.8735&1.3783&0.3890&1.4990&1.1122&1.3093&1.6269&1.0214\\
Ex2 & 1.0709&1.3432&0.5605&1.4621&1.4575&1.2719&2.4228&1.0334\\
Ex3 &0.9950&1.3793&0.5010&1.4779&1.3538&1.2747&2.2439&1.0334\\
\hline
\end{tabular}}
\caption{Rate of growth}
\end{table}

\begin{table}[htb] {\centering{\small
\begin{tabular}{|c|r|rrrr|}\hline
&\multicolumn{1}{|c}{n} &\multicolumn{1}{|c} {Alg 1} & \multicolumn{1}{c}{Alg 2} & \multicolumn{2}{c|}{Alg 3}\\ \hline
&400 &222.0 & 336.4 & 158.2 & 20.8\\
&800 & 386.6 & 481.0 & 185.8 & 24.4\\
Ex1&1600 &
502.0 & 597.2 & 210.8 & 27.6\\
&3200 &
644.6 & 670.2 & 225.4 & 28.6\\
&6400 &
745.2 & 781.4 & 244.2 & 30.2\\ \hline
&400& 302.4& 345.0 & 247.2 & 27.6\\
&800& 456.6 & 536.0 & 294.6 &32.2\\
Ex2&1600 & 570.6 & 602.8 & 312.2 & 32.8\\
&3200 & 735.0& 713.4 & 330.0 & 34.6\\
&6400 & 869.6 & 794.2 & 360.0 & 38.0\\ \hline
&400& 297.4 & 334.0 & 224.6&17.8 \\
&800 &
449.8 & 491.8 & 261.0 & 19.8 \\
Ex3&1600&
553.8 & 545.6& 279.2& 21.0\\
&3200& 732.6 & 657.8 & 293.6 &20.6 \\
&6400& 838.6 & 699.4 & 304.2 & 20.8 \\ \hline
\end{tabular}}
 
\caption{
Number of steps in dependence of the number of inequalities}

}
\end{table}

 


In the tables, the second column for the third algorithm gives the number of rescalings. Moreover the numbers for the simplex algorithm do not give the $d$ steps necessary to obtain the starting tableau. We see that all algorithms show little dependence on the number of constraints, in particular for Alg 3 this dependence appears to be sub--logarithmic.

For the dependence on $d$ we may not expect spectacular results, as we know by the remarks at the beginning, that we may not hope for less than $d$ iterations and on the other hand the simplex algorithm solved all examples in less than $12d$ iterations. But we see that for large dimensions the new algorithms need fewer iterations and in particular Alg 3, the algorithm with the rescaling, did very well.

To obtain some more insight in the dependence on the dimension we did a best fit of the form $\alpha d^{\beta}$ for all series and all algorithms. Here we computed two values for the simplex algorithm, for the first value we counted the $d$ steps for the initialization while in the second value these steps were neglected. For the algorithm with the rescaling we did not take into account the numbers of rescaling, as this number appears to have little dependence on the dimension. The results are presented in the second table.

We find that the simplex algorithm and Algorithm \ref{algo2} are superlinear, where the exponent of Algorithm \ref{algo2} is slightly smaller for all series, while the experimental algorithm is very close to being linear for our examples and thus for these examples very close to the hypothetical Gauss--type algorithm. 

Altogether the examples show that the concepts developed here might be of practical interest.

Finally the author thanks {\sc M. Henk} and {\sc A. Schürmann} for valuable discussions, hints to the literature and reading a previous version of this paper.

\end{document}